\documentclass[11pt]{amsart}
\usepackage[T1]{fontenc}
\usepackage[utf8]{inputenc}
\usepackage{lmodern}
\usepackage[margin=1in]{geometry}
\usepackage{amsmath,amssymb,mathtools,mathrsfs}
\usepackage{microtype,enumitem}
\usepackage[colorlinks=true,linkcolor=blue,citecolor=blue,urlcolor=blue]{hyperref}
\usepackage[capitalise,noabbrev]{cleveref}
\usepackage{tikz}
\usetikzlibrary{arrows.meta,positioning}

\newtheorem{theorem}{Theorem}[section]
\newtheorem{lemma}[theorem]{Lemma}
\newtheorem{proposition}[theorem]{Proposition}
\newtheorem{corollary}[theorem]{Corollary}
\theoremstyle{definition}
\newtheorem{definition}[theorem]{Definition}
\newtheorem{question}[theorem]{Question}

\theoremstyle{remark}
\newtheorem{remark}[theorem]{Remark}
\usepackage[capitalise,noabbrev]{cleveref}

\crefname{lemma}{lemma}{lemmas}
\Crefname{lemma}{Lemma}{Lemmas}

\crefname{proposition}{proposition}{propositions}
\Crefname{proposition}{Proposition}{Propositions}

\crefname{corollary}{corollary}{corollaries}
\Crefname{corollary}{Corollary}{Corollaries}

\crefname{definition}{definition}{definitions}
\Crefname{definition}{Definition}{Definitions}

\crefname{question}{question}{questions}
\Crefname{question}{Question}{Questions}

\crefname{example}{example}{examples}
\Crefname{example}{Example}{Examples}

\crefname{remark}{remark}{remarks}
\Crefname{remark}{Remark}{Remarks}
\newcommand{\silp}{\operatorname{silp}}
\newcommand{\spli}{\operatorname{spli}}
\newcommand{\pd}{\operatorname{pd}}
\newcommand{\idim}{\operatorname{id}}
\newcommand{\FPD}{\operatorname{FPD}}
\newcommand{\FID}{\operatorname{FID}}

\newcommand{\Ext}{\operatorname{Ext}}
\newcommand{\Hom}{\operatorname{Hom}}

\newcommand{\op}{\mathrm{op}}

\title[Rings with unequal silp--spli]{Rings with unequal \texorpdfstring{$\silp$--$\spli$}{silp--spli}}
\author{Xuyang Chen}
\email{xuyangchen19@fudan.edu.cn}
\author{Honghui Ding}
\email{2021000768@ruc.edu.cn}
\date{\today}

\begin{document}
\begin{abstract}
For a unital associative ring $R$, let $\silp(R)$ be the supremum of the injective dimensions of all projective right $R$-modules, and let $\spli(R)$ be the supremum of the projective dimensions of all injective right $R$-modules. We construct an explicit non-Noetherian unital algebra \[R=k1\oplus
\frac{k(0\longrightarrow1\longrightarrow2\longrightarrow\cdots)}
{\langle\text{all paths of length }2\rangle}\] over any field and prove that
\(
\silp(R)=1\) and \(\spli(R)=\infty.
\)
The algebra is the unitization of the non-unital path algebra of the one-sided infinite ray modulo paths of length two. The example answers the question whether $\silp=\spli$ raised by Gedrich and Gruenberg.
\end{abstract}
\maketitle

\section{Introduction}
For a unital associative ring $R$, a right $R$-module is \emph{projective} if it is a direct summand of a free right module, and it is \emph{injective} if every homomorphism into it defined on a submodule extends to the ambient module. A projective resolution of a right module $M$ is an exact sequence of projective modules ending in a surjection onto $M$; its smallest possible length, or infinity if no finite length is possible, is the \emph{projective dimension} $\pd_R M$. The \emph{injective dimension} $\idim_R M$ is defined dually. We set
\begin{equation}\label{eq:invariants}
\begin{split}
\silp(R)&=\sup\{\idim_R P:\ P\text{ projective right }R\text{-module}\},\\
\spli(R)&=\sup\{\pd_R I:\ I\text{ injective right }R\text{-module}\}.
\end{split}
\end{equation}
All dimensions are taken in $\mathbb{N}\cup\{\infty\}$, with $0\in\mathbb{N}$. We make the right-module convention explicit throughout. 

Gedrich and Gruenberg introduced these invariants in their study of complete cohomological functors on groups \cite{GG1987}. Their Theorems~1.5 and~$1.5'$ describe terminal completions of Ext in terms of the respective finiteness conditions; immediately following these results, on page~209, they ask whether one invariant can be finite while the other is infinite. Their Proposition~1.6 proves that the invariants agree if both are finite, and Corollary~1.7 identifies both with the global dimension when it is finite \cite[pp.~209--210]{GG1987}. The authors also recall that either invariant is zero precisely for quasi-Frobenius rings, and discuss positive implications for commutative Noetherian rings and certain Hopf and group algebras \cite[pp.~203--204, \S2]{GG1987}.

Auslander and Reiten investigated Gorenstein Artin algebras in \cite{AR1991}. In their setting, two-sided finiteness of the self-injective dimension of the regular module defines the Gorenstein property; see especially \cite[p.~223]{AR1991}. The \emph{Gorenstein symmetry conjecture} asks whether finiteness on one side forces finiteness on the other. For Artin algebras, this is equivalent to the one-sided $\silp$--$\spli$ question, although the two questions have different scopes outside that class; see \cite[p.~408]{ET2011} and \cite[Chapter~VII, p.~123]{BR2007}. Beligiannis developed a broad homological theory of Gorenstein categories \cite{Bel2000}; Beligiannis and Reiten explicitly asked whether either finite invariant alone makes an abelian category with enough projectives and injectives Gorenstein \cite[Chapter~VII, Definition~2.1 and the note after Remark~2.7]{BR2007}.

Emmanouil and Talelli returned to the ring-theoretic question in \cite[pp.~408--410]{ET2011}. They introduced comparisons using the supremum of the flat dimensions of injective modules, and distinguished left and right invariants carefully. Dalezios and Emmanouil subsequently obtained comparisons between invariants of $R$ and $R^{\op}$ via a class of Gorenstein flat modules \cite{DE2023}. Christensen, Estrada, Liang, Thompson and Wang established more general identities linking the Gorenstein global dimension with $\silp$, $\spli$ and the large finitistic dimensions \cite[Theorem~1.18]{CELTW2024}. Chatzistavridis and Emmanouil proved a comparison formula for two-sided weakly coherent rings \cite[Theorem~3.4]{CE2025}.

Wang, Li, Hu and Zhu explicitly describes the equality for arbitrary unital rings as an open question in \cite[p.~3]{WLHZ2026}. Our argument below is independent of their results.

\begin{question}[Gedrich--Gruenberg \cite{GG1987}]\label{q:GG}
Can one have $\silp(R)<\infty$ and $\spli(R)=\infty$, or $\spli(R)<\infty$ and $\silp(R)=\infty$, for a unital associative ring $R$?
\end{question}

We answer the first alternative affirmatively, which is sufficient to disprove a universal equality. 

\begin{theorem}\label{thm:main}
Let $k$ be an arbitrary field. There exists an associative unital $k$-algebra $R$, with a countable basis over $k$, such that
\[
\silp(R)=1\qquad\text{and}\qquad\spli(R)=\infty.
\]
Thus the extended-integer equality $\silp(R)=\spli(R)$ fails for unital rings in general.
\end{theorem}

The proof has three essential parts. First, the modules generated by a certain augmented ideal $A$ of $R$ are identified with nonnegative chain complexes of vector spaces and admit explicit decompositions. Second, a canonical short exact sequence for any right module gives a uniform bound on the injective dimension of every free module. Finally, an injective simple has infinite projective dimension. We give every construction and all dimension calculations below.

\section{The algebra and its basic modules}\label{sec:algebra}
Fix a field $k$. Let $Q$ be the one-sided infinite ray
\[
Q:\qquad
\begin{tikzpicture}[baseline=(current bounding box.center),>=Latex,node distance=18mm]
\node (0) {$0$};
\node (1) [right=of 0] {$1$};
\node (2) [right=of 1] {$2$};
\node (3) [right=of 2] {$3$};
\node (4) [right=of 3] {$\cdots$};
\draw[->] (0) -- node[above] {$a_0$} (1);
\draw[->] (1) -- node[above] {$a_1$} (2);
\draw[->] (2) -- node[above] {$a_2$} (3);
\draw[->] (3) -- (4);
\end{tikzpicture}
\]
and let
\[
A=kQ/\langle\text{all paths of length }2\rangle.
\]
Equivalently, $A$ is a rng (or a non-unital ring) whose underlying $k$-vector space has basis $e_i,a_i$ ($i\geq0$), where $e_i$ is the trivial path at vertex $i$ and $a_i$ is the arrow from $i$ to $i+1$, with multiplication
\begin{equation}\label{eq:mult}
 e_ie_j=\delta_{ij}e_i,\qquad
 e_ia_j=\delta_{ij}a_j,\qquad
 a_je_i=\delta_{i,j+1}a_j,\qquad
 a_ia_j=0.
\end{equation}
The relations express path concatenation, with all length-two paths declared zero; in particular multiplication is associative. Every element of $A$ is a \emph{finite} linear combination of the displayed basis elements. The algebra $A$ has no identity. Its \emph{unitization} (or Dorroh extension) is the unital algebra
\begin{equation}\label{eq:unitize}
R=k1\oplus A
= k1\oplus \frac{kQ}{\langle\text{all paths of length }2\rangle},
\qquad
(\lambda1+x)(\mu1+y)=\lambda\mu1+\lambda y+\mu x+xy.
\end{equation}
The element $1$ is the identity, $A$ is a two-sided ideal, and $R/A\cong k$ as rings.

Since $e_i=e_i^2$ and $a_i=e_ia_i$, we have $A^2=A$. Moreover, as a right module,
\begin{equation}\label{eq:Adec}
A=\bigoplus_{i\geq0}P_i,\qquad
P_i=e_iR=ke_i\oplus ka_i.
\end{equation}
Indeed, the summands are independent, have union spanning $A$, and all sums have finite support. Each $P_i$ is projective because $e_i$ is an idempotent and $e_iR$ is a direct summand of $R_R$. Thus $A_R$ is projective.

Define the one-dimensional right module $S_i=ks_i$ by
\[
s_ie_j=\delta_{ij}s_i,\qquad s_ia_j=0,\qquad s_i1=s_i.
\]
The map $P_i\to S_i$ defined by $e_i\mapsto s_i$ and $a_i\mapsto0$ is surjective and has kernel $ka_i\cong S_{i+1}$, giving
\begin{equation}\label{eq:ses}
0\longrightarrow S_{i+1}
\longrightarrow P_i
\longrightarrow S_i\longrightarrow0.
\end{equation}
The left-hand inclusion sends $s_{i+1}$ to $a_i$. Finally, set $S_\infty=R/A$, the one-dimensional module on which $A$ acts as zero. From \eqref{eq:Adec},
\begin{equation}\label{eq:Sinf}
0\longrightarrow A\longrightarrow R\longrightarrow S_\infty\longrightarrow0
\end{equation}
is a projective resolution of length at most one, so
\begin{equation}\label{eq:pdSinf}
\pd_R S_\infty\leq1.
\end{equation}

\section{Modules generated by the vertex ideal}\label{sec:unitary}
\begin{definition}
A right $R$-module $N$ is called \emph{$A$-unitary} if $NA=N$. It need not be finitely generated.
\end{definition}

For a module $M$ and a cardinal $\alpha$, we denote by $M^{(\alpha)}$ the direct sum of a set $I$ of copies of $M$ with $|I|=\alpha$.

\begin{lemma}\label{lem:complex}
If $N$ is an $A$-unitary right $R$-module, then there exist cardinals $\alpha_i,\beta_i$ such that
\begin{equation}\label{eq:split}
N\cong\bigoplus_{i\geq0}P_i^{(\alpha_i)}
\oplus\bigoplus_{i\geq0}S_i^{(\beta_i)}.
\end{equation}
The direct sums and cardinals may be infinite.
\end{lemma}
\begin{proof}
Every $n\in N=NA$ is a finite sum $\sum_t n_tb_t$ with $b_t\in A$. Every $b_t$ is supported at finitely many vertices and satisfies $b_t=\sum_i b_te_i$. Consequently $n=\sum_i ne_i$ with only finitely many nonzero terms. The orthogonality of the $e_i$ then implies an internal direct sum of $k$-vector spaces
\[
N=\bigoplus_{i\geq0}V_i,\qquad V_i=Ne_i.
\]
These $V_i$ need not individually be $R$-submodules. Right multiplication by $a_i$ defines a $k$-linear map $d_i:V_i\to V_{i+1}$. The relations $a_i a_{i+1}=0$ imply $d_{i+1}d_i=0$, so the sequence
\[
V_0\xrightarrow{d_0}V_1\xrightarrow{d_1}V_2\longrightarrow\cdots
\]
is a complex of vector spaces.

Put $B_0=0$ and $B_i=\operatorname{im}d_{i-1}$ for $i\geq1$. For each $i$, choose vector-space complements $H_i$ and $C_i$ such that
\[
\ker d_i=B_i\oplus H_i,
\qquad V_i=B_i\oplus H_i\oplus C_i.
\]
The restriction $d_i\vert_{C_i}$ maps isomorphically onto $B_{i+1}$. Choose bases of all $H_i$ and $C_i$ and transport the corresponding bases to all $B_i$ through these isomorphisms. Each basis vector of $H_i$ spans a submodule isomorphic to $S_i$, while a vector in $C_i$ together with its image in $B_{i+1}$ spans a submodule isomorphic to $P_i$. These submodules are independent, stable under the actions in \eqref{eq:mult}, and together exhaust $N$. The resulting direct-sum decomposition is \eqref{eq:split}.
\end{proof}

\begin{lemma}\label{lem:Extvanish}
If $F=R^{(\kappa)}$ is a free right module of arbitrary rank $\kappa$, and $N$ is $A$-unitary, then
\begin{equation}\label{eq:Extvanish}
\Ext_R^q(N,F)=0\qquad\text{for every }q\geq1.
\end{equation}
\end{lemma}
\begin{proof}
It is enough by \cref{lem:complex} to prove the vanishing for $S_i$, since $P_i$ is projective. Indeed, $\Ext^q_R(\oplus_{i\in I}N_i,F)\cong\prod_{i\in I}\Ext_R^q(N_i,F)$.

Splicing \eqref{eq:ses} gives a projective resolution of $S_i$:
\begin{equation}\label{eq:res}
\cdots\longrightarrow P_{i+2}\longrightarrow P_{i+1}
\longrightarrow P_i\longrightarrow S_i\longrightarrow0.
\end{equation}
For $j\geq1$, the differential $P_j\to P_{j-1}$ takes $e_j$ to $a_{j-1}$ and $a_j$ to zero; its image is the kernel of the map $P_{j-1}\to P_{j-2}$ (or, at the last step, of $P_i\to S_i$). Thus \eqref{eq:res} is exact.

Apply $\Hom_R(-,F)$ and we get a cochain complex\[
0\longrightarrow \Hom(P_i,F)\longrightarrow\Hom_R(P_{i+1},F)\longrightarrow\Hom_R(P_{i+2},F)\longrightarrow\cdots.
\]
For each $j$, evaluation at $e_j$ gives a natural isomorphism of right $e_jRe_j$-modules $\Hom_R(P_j,F)=\Hom_R(e_jR,F)\cong Fe_j$. The cochain differential induced by $P_{j+1}\to P_j$ becomes right multiplication by $a_j$,
\[
Fe_j\xrightarrow{\,\cdot a_j\,}Fe_{j+1}.
\]
Since $Re_j=ke_j\oplus ka_{j-1}$ for $j\geq1$, we have
\[
Fe_j=(ke_j\oplus ka_{j-1})^{(\kappa)}\qquad(j\geq1).
\]
The rules $e_ja_j=a_j$ and $a_{j-1}a_j=0$ therefore show
\begin{equation}\label{eq:kernelimage}
\ker(Fe_j\xrightarrow{\cdot a_j}Fe_{j+1})
=(ka_{j-1})^{(\kappa)}
=\operatorname{im}(Fe_{j-1}\xrightarrow{\cdot a_{j-1}}Fe_j)
\end{equation}
for every $j\geq1$. By \eqref{eq:kernelimage}, all those cohomology groups of degree above 1 vanish, and hence $\Ext_R^q(S_i,F)=0$ for $q\geq1$. The direct-sum decomposition \eqref{eq:split} finishes the proof.
\end{proof}

\begin{lemma}\label{lem:filtration}
Every right $R$-module $M$ has a short exact sequence
\begin{equation}\label{eq:filtration}
0\longrightarrow MA\longrightarrow M\longrightarrow S_\infty^{(\lambda)}\longrightarrow0
\end{equation}
for some cardinal $\lambda$, and $MA$ is $A$-unitary.
\end{lemma}
\begin{proof}
Since $A$ is a two-sided ideal, $MA$ is a right submodule. The identity $A^2=A$ implies $(MA)A=MA$. The quotient $M/MA$ is annihilated by $A$, so its action factors through $R/A\cong k$. Choose a vector-space basis of this quotient; each basis vector generates a copy of $S_\infty$. This gives \eqref{eq:filtration}.
\end{proof}

\begin{proposition}\label{prop:upper}
Every projective right $R$-module has injective dimension at most one.
\end{proposition}
\begin{proof}
Let $F=R^{(\kappa)}$ be free, of any rank $\kappa$, and let $M$ be an arbitrary right module. Apply $\Hom_R(-,F)$ to \eqref{eq:filtration} and we get a long exact sequence\[
\dots\longrightarrow\Ext_R^i(S_\infty^{(\lambda)},F)\longrightarrow \Ext_R^i(M,F)\longrightarrow \Ext_R^i(MA,F)\longrightarrow \Ext_R^{i+1}(S_\infty^{(\lambda)},F)\longrightarrow\dots.
\] By \cref{lem:Extvanish}, $\Ext_R^q(MA,F)=0$ for all $q\geq1$. Since $\pd_R S_\infty\leq1$ by \eqref{eq:pdSinf}, any direct sum of copies of $S_\infty$ has projective dimension at most one. Therefore $\Ext_R^q(S_\infty^{(\lambda)},F)=0$ for $q\geq2$. The long exact sequence gives
\[
\Ext_R^q(M,F)=0\qquad\text{for every }M\text{ and every }q\geq2.
\]
For a module $F$, vanishing of these groups for all first arguments is equivalent to $\idim_R F\leq1$. Finally, every projective right module is a direct summand of some free module, and injective dimension cannot increase under direct summands. The assertion follows.
\end{proof}

\section{The injective simple of infinite projective dimension}\label{sec:infinite}
\begin{lemma}\label{lem:injective}
The simple right module $S_0$ is injective.
\end{lemma}
\begin{proof}
The left ideal $Re_0$ has basis $e_0$, because no arrow ends at vertex $0$. It is projective as a left $R$-module, since $e_0$ is an idempotent. For a left $R$-module $L$, let $D(L)=\Hom_k(L,k)$ with the usual right action $(\varphi\cdot r)(x)=\varphi(rx)$. For any right module $M$ we have natural isomorphisms of vector spaces
\[
\Hom_R(M,D(Re_0))\cong D(M\otimes_R Re_0)\cong D(Me_0).
\]
The functor to the vector spaces $M\mapsto Me_0$ is exact: on an exact sequence of right modules, multiplication by the idempotent $e_0$ preserves kernels and surjections. Duality $D$ is exact on vector spaces. Hence $\Hom_R(-,D(Re_0))$ is exact, and $D(Re_0)$ is injective. It is one-dimensional, $e_0$ acts as the identity, and all arrows and other vertex idempotents act as zero; thus $D(Re_0)\cong S_0$.
\end{proof}

\begin{proposition}\label{prop:lower}
The injective right module $S_0$ has infinite projective dimension. In particular, $\spli(R)=\infty$.
\end{proposition}
\begin{proof}
Use the projective resolution \eqref{eq:res} with $i=0$. The term in homological degree $j$ is $P_j$. For $n\geq1$, the evaluation isomorphism gives
\[
\Hom_R(P_j,S_n)\cong S_ne_j=
\begin{cases}k,&j=n,\\0,&j\ne n.\end{cases}
\]
Thus $\Hom_R(P_\bullet,S_n)$ has exactly one nonzero cochain group, equal to $k$ in degree $n$. It follows that
\begin{equation}\label{eq:nonzero}
\Ext_R^n(S_0,S_n)\cong k\ne0\qquad(n\geq1).
\end{equation}
Since \eqref{eq:nonzero} holds for unbounded $n$, $\pd_R S_0=\infty$. By \cref{lem:injective}, $S_0$ is injective, and the assertion follows from \eqref{eq:invariants}.
\end{proof}

\section{Sharp bound for injective dimension}\label{sec:sharp}
\begin{proposition}\label{prop:sharp}
The regular right module $R_R$ has injective dimension exactly one. Consequently $\silp(R)=1$.
\end{proposition}
\begin{proof}
By \cref{prop:upper}, $\idim_R R\leq1$. To see that $R_R$ is not injective, define $c_i=0$ for even $i$ and $c_i=1$ for odd $i$. Since $A=\bigoplus_i P_i$, there is a well-defined right $R$-linear map $f:A\to R$ given on each $P_i$ by
\[
f(e_i r)=c_i e_i r\qquad(r\in R).
\]
Suppose that $f$ extends to an endomorphism $\widetilde f:R_R\to R_R$. Such an endomorphism is left multiplication by $r_0=\widetilde f(1)$. Write $r_0=\lambda1+b$ with $b\in A$. Since $b$ has finite support, $be_i=0$ for all sufficiently large $i$. Therefore, for all sufficiently large $i$,
\[
f(e_i)=\widetilde f(e_i)=r_0e_i=\lambda e_i.
\]
This is impossible because $f(e_i)=c_ie_i$ and $(c_i)$ is not eventually constant. Thus $f$ cannot be extended, so $R_R$ is not injective and $\idim_R R\geq1$. Together with the upper bound this gives $\idim_R R=1$. Since $R_R$ is projective, $\silp(R)\geq1$; \cref{prop:upper} gives the opposite inequality.
\end{proof}

\begin{proof}[Proof of \cref{thm:main}]
The construction \eqref{eq:mult}--\eqref{eq:unitize} gives an associative unital algebra with basis $1,e_0,a_0,e_1,a_1,\ldots$. By \cref{prop:sharp}, $\silp(R)=1$, and by \cref{prop:lower}, $\spli(R)=\infty$.
\end{proof}

\begin{corollary}\label{cor:fid}
For this ring, the big finitistic projective and injective dimensions satisfy
\[
\FPD(R)=1\qquad\text{and}\qquad\FID(R)=\infty.
\]
Here $\FPD$ (respectively $\FID$) is the supremum of finite projective (respectively injective) dimensions attained by arbitrary right modules.
\end{corollary}
\begin{proof}
If a right module $X$ has finite projective dimension $d>0$, the last short exact sequence in a minimal-length projective resolution does not split. Thus $\Ext_R^d(X,P)\ne0$ for some projective $P$ (the leftmost projective in that resolution). Consequently $d\leq\silp(R)=1$. This proves $\FPD(R)\leq1$. On the other hand, $\pd_R S_\infty\leq1$ by \eqref{eq:pdSinf}, and it is not zero: if \eqref{eq:Sinf} split, every homomorphism $A\to R$ would extend to an endomorphism of $R$, contradicting the map used in \cref{prop:sharp}. Hence $\pd_R S_\infty=1$ and $\FPD(R)=1$.

For the second assertion, put $I_0=S_0$ and $I_i=D(Re_i)$ for $i\geq1$. The same argument as in \cref{lem:injective} shows each $I_i$ injective. Directly from \eqref{eq:mult}, $I_i\cong P_{i-1}$ for $i\geq1$. Therefore \eqref{eq:ses} gives $0\to S_i\to I_i\to S_{i-1}\to0$, and induction on $i$, starting with injective $S_0$, yields $\idim_R S_i\leq i$. Conversely \eqref{eq:nonzero} gives $\idim_R S_i\geq i$. Thus $\idim_R S_i=i$ for every $i\geq0$, proving $\FID(R)=\infty$.
\end{proof}

\begin{remark}
The four values $\silp(R)=\FPD(R)=1$ and $\spli(R)=\FID(R)=\infty$ are consistent with the max formulas of Christensen--Estrada--Liang--Thompson--Wang \cite[Theorem~1.18]{CELTW2024}.
\end{remark}

\begin{remark}[$R$ is non-Noetherian]
The right ideals $\bigoplus_{i=0}^m e_iR$ form a strictly increasing chain, so $R$ is not right Noetherian. Similarly, the left ideals $\bigoplus_{i=0}^m Re_i$ form a strictly increasing chain, so it is not left Noetherian.
\end{remark}

\end{document}